\documentclass[12pt]{amsart}
\usepackage[english]{babel}
\usepackage{amsfonts,amssymb,latexsym,amscd}

\theoremstyle{plain}
\newtheorem{theorem}{Theorem}

\newtheorem{proposition}{Proposition}

\theoremstyle{definition}
\newtheorem{definition}{Definition}
\theoremstyle{remark}
\newtheorem{remark}{Remark}

\usepackage{hyperref}

\begin{document}

\title{Movable categories}
\author{Pavel S. Gevorgyan}
\address{Moscow Pedagogical State University}
\email{pgev@yandex.ru}

\begin{abstract}
The notion of movability for metrizable compacta was introduced by K.Borsuk. In this paper we define the notion of a movable category and prove that the movability of a topological space $X$ coincides with the movability of a suitable category, which is generated by the topological space $X$ (i.e., the category $\mathcal{W}^X$, defined by S.Marde\v{s}i\'{c}). 
\end{abstract}



\maketitle

\section{Introduction}

The notion of movability for metrizable compacta was introduced by K. Borsuk \cite{borsuk}. For the more general cases this notion was
extended by S. Marde\v{s}i\'{c} and J. Segal \cite{mardsegal2}, J.
Segal~\cite{segal}, P. Shostak~\cite{shostak}. In the equivariant theory of
shape this notion was studied in the works of the author of
present article~$[4,5,6,7,8]$ and of Z. \v{C}erin~\cite{cerin}.

It is necessary to note that in the works mentioned above the
movability of topological spaces was defined by means of
neighborhoods of the given space (embedded as closed set in a
certain $AR$-space) or by means of inverse systems, depending on
the approach to shape theory used. However, the categorical
approach to shape theory  of S.~Marde\v{s}i\'{c}~\cite{mardsegal1} lacks a
suitable categorical definition of movability.

In this article we define the notion of movable category and prove
that the movability of a topological space is equivalent to the
movability of a certain category.

The author is extremely grateful to the referee for his helpful
remarks and comments.

\section{Basic notions and conventions concerning shape and movability.}

Let $\mathcal{HTOP}$ denote the homotopy category of topological
spaces and homotopy classes of maps and $\mathcal{HCW}$ the full
subcategory of $\mathcal{HTOP}$ whose objects are all topological
spaces having the homotopy type of a CW-complex.

\begin{definition}[{K. Morita [11]}] An inverse system
$\{X_\alpha , p_{\alpha \alpha '}, A\}$ in $\mathcal{HCW}$ is
called {\bf associated} with or an {\bf expansion} of a
topological space $X$ if there are homotopy classes $p_\alpha : X
\to X_\alpha $ for $\alpha \in A$ such that the following
conditions are satisfied.

(1) \ $p_{\alpha \alpha'}p_{\alpha'}=p_\alpha$, if
$\alpha<\alpha'$.

(2) \ For any homotopy class $f:X\to Q$ with $Q\in$
Ob($\mathcal{HCW}$), there exists $\alpha \in A$ and a homotopy
class $f_\alpha:X_\alpha \to Q$ such that $f=f_\alpha p_\alpha$.

(3) \ For $\alpha \in A$ and for homotopy classes $f_\alpha,
g_\alpha :X_\alpha \to Q$ with $Q\in$ Ob($\mathcal{HCW}$) such
that $f_\alpha p_\alpha=g_\alpha p_\alpha$, there exist $\alpha '
\in A$ with $\alpha \leqslant \alpha '$ such that $f_\alpha
p_{\alpha \alpha '}=g_\alpha p_{\alpha \alpha '}$.

\end{definition}

For any topological space $X$ there exist an inverse system in
$\mathcal{HCW}$ associated with $X$ [11], and so there is a shape
theory for all topological spaces because the abstract theory of
shape yields that there is a shape theory for $\mathcal{HTOP}$ iff
every topological space has an expansion, i.e., iff
$\mathcal{HCW}$ is a so-called "dense" subcategory of
$\mathcal{HTOP}$.

On the categorical approach to shape theory of S.
Marde\v{s}i\'{c}~\cite{mardsegal1} for each topological space $X$ it is
introduced a new comma category $\mathcal{W}^X$, whose objects are
homotopy classes $f : X \to Q$ and whose morphisms are the
following commutative triangles \\

\begin{center}
\begin{picture}(80,40)
\put(0,0){$Q'$} \put(40,40){$X$} \put(80,0){$Q$}
\put(38,38){\vector(-1,-1){28}} \put(52,38){\vector(1,-1){28}}
\put(74,3){\vector(-1,0){60}} \put(14,21){$f'$} \put(73,21){$f$}
\put(40,-4){$\mu $}
\end{picture}
\end{center}
\quad \\ where $Q, Q' \in$ Ob($\mathcal{HCW}$). Then a shape map
$f:X\to Y$ is defined as a covariant functor $f:\mathcal{W}^Y \to
\mathcal{W}^X $ which keeps morphisms $\mu$ fixed.

\begin{definition}[{[10]}]
An inverse system $\{X_\alpha , p_{\alpha \alpha '}, A\}$ in
$\mathcal{HCW}$ is called movable if

(*) for every $\alpha \in A$, there exists an $\alpha ' \in A$,
$\alpha ' \geqslant \alpha $ such that for all $\alpha '' \in A$,
$\alpha '' \geqslant \alpha $, there exists a homotopy class
$r^{\alpha ' \alpha ''}:X_{\alpha '}\to X_{\alpha ''}$ such that
$$p_{\alpha \alpha'}=p_{\alpha \alpha''}r^{\alpha ' \alpha ''} .$$

The topological space $X$ is called movable if there exist an
inverse system $\{X_\alpha , p_{\alpha \alpha '}, A\}$ in
$\mathcal{HCW}$ which associated with $X$ and which is movable.
\end{definition}

The reader is referred to the book by S. Marde\v{s}i\'{c} and J.
Segal [9] for general information about shape theory.

\section{The movable categories.}

Let $K$ be an arbitrary category and $K'$ any subcategory of the
category $K$.
\begin{definition}\label{movcat}
We say that a subcategory $K'$ is movable in a category $K$, if
for any object $X\in$ Ob($K'$)   there exists an object $Y\in$
Ob($K'$) and a morphism $f~\in~ K'(Y, X)$  such that for any
object $Z\in$ Ob($K'$)  and any morphism $g\in K'(Z, X)$  there is
a morphism $h\in K(Y, Z)$  which make the following diagram
commutative \\[10pt]

\end{definition}

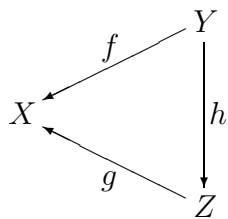
\begin{center}
\begin{picture}(70,70)
\put(0,35){$X$}
\put(70,0){$Z$}
\put(70,70){$Y$}

\put(67,71){\vector(-2,-1){54}}
\put(35,61){$f$}
\put(67,7){\vector(-2,1){54}}
\put(35,12){$g$}
\put(74,66){\vector(0,-1){55}}
\put(76,35){$h$}
\end{picture}
\end{center}

\quad

\begin{definition}

 We say that a category is movable if it is movable in  itself.

\end{definition}

\begin{definition}[{\cite{bukur}}]\label{zeromorp}
 It is said that $K$ is a category with zero-morphisms if for any pair $(A, B)$
 of objects from a category $K$ there exist morphisms $o_{BA} : A\to B$
 which, for all morphisms $\nu : B \to C$ and $u:D \to A$, where $C$ and $D$
 are objects of the category $K$, satisfy the following equalities
\[
\nu o_{BA}=o_{CA}, \qquad o_{BA} u=o_{BD}.
\]
\end{definition}

\begin{definition}[{\cite{bukur}}]

 An object $O\in$ Ob($K$) is called initial if for any object $X\in$ Ob($K$)
 the set $Mor_K(O, X)$ consists of a single morphism.

\end{definition}

\begin{proposition}
 Any category $K$ with zero-morphisms is movable.
\end{proposition}

\begin{proof}
 Let $X\in$ Ob($K$)  be an arbitrary object. It appears that for
 the object we seek (see definition \ref{movcat}), we may take any object $Y\in ob(K)$
 and for the morphism $f\in K(Y, X)$  it is necessary to take a zero-morphism
 $o_{XY} : Y\to X$. Indeed, let $g\in K(Z, X)$ be an arbitrary morphism.
 It is clear that zero-morphism $o_{ZY} : Y\to~Z$ is the morphism we seek,
 that is $go_{ZY}=o_{XY}$, which follows from definition \ref{zeromorp}.
\end{proof}

\begin{proposition}
 Any category $K$ with initial objects is movable.
\end{proposition}

\begin{proof}
 Let $X\in$ Ob($K$) be any object. Let us consider the initial object $O$
 of the category $K$. We denote by $u_X$ the single morphism from the object $O$ to
 the object $X$. Now it is not difficult to note that the object $O$ and
 the morphism $u_X : O \to X$ satisfies the condition of definition \ref{movcat}.
 Indeed: let $Y\in ob(K)$ be any object and $g : Y \to X$ be any morphism of the category $K$.
 It is clear that the single morphism $u_Y : O \to Y$ satisfies the condition $u_X=g\circ u_Y$.
\end{proof}

\section{The movability of topological spaces.}

\begin{theorem}
 The topological space $X$ is movable if and only if the category $\mathcal{W}^X$ is
 movable.
 \end{theorem}

This theorem is a simple reformulation of the following theorem.

\begin{theorem}
 The topological space $X$ is movable if and only if the following condition
 is satisfied.

\hangindent=0.5cm (*) For any $Q\in$ Ob($\mathcal{HCW}$) and any
homotopy class $f:X\to Q$ there exist
 $Q'\in$ Ob($\mathcal{HCW}$) and homotopy classes
 $f':X\to Q'$, $\eta :Q'\to Q$,
 satisfying $f=\eta \circ f'$, such that for any $Q''\in$ Ob($\mathcal{HCW}$) and
 homotopy classes $f'':X\to Q''$, $\eta ' :Q''\to Q$, satisfying the
 condition  $f=\eta '\circ f''$, there exist a homotopy class
 $\eta '':Q'\to Q''$ which satisfies the condition $\eta =\eta '\circ \eta ''$ (diagram
 1).
\end{theorem}
\begin{center}
\begin{picture}(60,100)
\put(0,0){$Q'$}
\put(0,60){$Q$}
\put(50,30){$X$}
\put(100,30){$Q''$}
\put(4,12){\vector(0,1){40}}
\put(13,4){\vector(4,1){81}}
\put(94,38){\vector(-3,1){80}}
\put(46,28){\vector(-2,-1){32}}
\put(45,39){\vector(-2,1){32}}
\put(59,33){\vector(1,0){33}}
\put(-4,30){$\eta$}
\put(24,22){$f'$}
\put(24,38){$f$}
\put(58,53){$\eta '$}
\put(67,36){$f''$}
\put(58,8){$\eta ''$}
\end{picture}

$ $ \linebreak
Diagram 1.

\end{center}

\begin{proof}
 Let condition $(*)$ be satisfied. We must prove that $X$ is movable.
 Let us consider an inverse system $\{X_\alpha , p_{\alpha \alpha '}, A\}$ in $\mathcal{HCW}$
 which is associated with the topological space $X$.

 Let $\alpha \in A$ be any element and $p_\alpha :X\to X_\alpha $ be the natural projection.
 By $(*)$ for the homotopy class
 $p_\alpha :X\to X_\alpha $ let $Q'\in$ Ob($\mathcal{HCW}$),
 and $f':X\to Q'$, $\eta :Q'\to X_\alpha $ are homotopy classes satisfying the condition
 $\eta \circ f'=p_\alpha $ (diagram~2).

 Since inverse system $\{X_\alpha , p_{\alpha \alpha '}, A\}$ is associated with $X$, there exists
 $\tilde{\alpha} \in A$, $\tilde{\alpha} \geqslant \alpha $
 and $\tilde{f'}:X_{\tilde {\alpha }} \to Q'$ such that
\begin{equation}
\label{1}
f'=\tilde{f'}\circ p_{\tilde{\alpha }}.
\end{equation}
It is not difficult to verify that
\begin{equation}
\label{2}
p_{\alpha \tilde {\alpha }} \circ p_{\tilde {\alpha }}=
\eta \circ \tilde{f'} \circ p_{\tilde{\alpha }}.
\end{equation}
Indeed:
\[
\eta \circ \tilde{f'} \circ p_{\tilde{\alpha }}=\eta \circ f'=p_\alpha=
p_{\alpha \tilde {\alpha }} \circ p_{\tilde {\alpha }}.
\]

From the equality (\ref{2}) and the definition 1 we infer the
existence of an index $\alpha ' \in A$, $\alpha '\geqslant
\tilde{\alpha}$ for which
\begin{equation}
\label{3}
p_{\alpha \tilde {\alpha }} \circ p_{\tilde{\alpha }\alpha '}=
\eta \circ \tilde{f'} \circ p_{\tilde{\alpha }\alpha '}.
\end{equation}

The obtained index $\alpha ' \in A$ satisfies the condition of the
movability of inverse system $\{X_\alpha , p_{\alpha \alpha '},
A\}$. Indeed, let $\alpha '' \in A$, $\alpha ''\geqslant \alpha$
be any element. For the homotopy classes $p_{\alpha \alpha
''}:X_{\alpha ''}\to X_\alpha$ and $p_{\alpha ''}:X\to X_{\alpha
''}$ (with the condition $p_\alpha=p_{\alpha \alpha ''}\circ
p_{\alpha ''}$) there exist a homotopy class $\eta '':Q' \to
X_{\alpha ''}$ , which satisfies the equality
\begin{equation}
\label{4}
\eta=p_{\alpha \alpha ''} \circ \eta ''.
\end{equation}
(see the condition $(*)$ ). Now it is easy to see that $g=\eta ''
\circ \tilde{f'} \circ p_{\tilde{\alpha} \alpha '}$ is the
homotopy class we seek, i. e. the following condition is
satisfied:
\begin{equation}
\label{5}
p_{\alpha \alpha '}=p_{\alpha \alpha ''} \circ g.
\end{equation}
Indeed:
\[
p_{\alpha \alpha '}=p_{\alpha \tilde{\alpha}} \circ
p_{\tilde{\alpha} \alpha '}=\eta \circ \tilde{f'} \circ
p_{\tilde{\alpha} \alpha '}=p_{\alpha \alpha ''} \circ \eta '' \circ
\tilde{f'} \circ p_{\tilde{\alpha} \alpha '}=p_{\alpha \alpha ''} \circ g.
\]

$ $\linebreak

\begin{center}
\begin{picture}(90,90)
\put(0,30){$X_\alpha$}
\put(0,60){$X_{\tilde{\alpha}}$}
\put(45,30){$X$}
\put(45,60){$X_{\alpha '}$}
\put(90,0){$X_{\alpha ''}$}
\put(90,90){$Q'$}
\put(87,6){\vector(-3,1){70}}
\put(54,30){\vector(3,-2){35}}
\put(43,34){\vector(-1,0){28}}
\put(55,38){\vector(2,3){32}}
\put(50,40){\vector(0,1){16}}
\put(5,56){\vector(0,-1){16}}
\put(43,63){\vector(-1,0){28}}
\put(87,88){\vector(-3,-2){30}}
\put(12,65){\vector(3,1){75}}
\put(94,84){\vector(0,-1){73}}

\put(-13,45){$p_{\alpha \tilde{\alpha}}$}
\put(36,45){$p_{\alpha '}$}
\put(95,45){$\eta ''$}
\put(20,57){$p_{\tilde{\alpha} \alpha '}$}
\put(22,37){$p_\alpha$}
\put(43,10){$p_{\alpha \alpha ''}$}
\put(45,80){$\tilde{f'}$}
\put(60,76){$\eta $}
\put(70,22){$p_{\alpha ''}$}
\end{picture} { } \\[15pt]

Diagram 2.
\end{center}

Now we must prove the converse. Let $X$ be a movable topological
space and some inverse system $\{X_\alpha , p_{\alpha \alpha '},
A\}$ associated with $X$ . Let us prove that the condition $(*)$
is satisfied. To this end, consider any homotopy class $f:X \to Q$
(diagram 3). From the association of the inverse system
$\{X_\alpha , p_{\alpha \alpha '}, A\}$ with the space $X$ follows
that there exist an index $\alpha \in A$ and a homotopy class
$f_\alpha :X_\alpha \to Q$ such that
\begin{equation}
\label{6}
f=f_\alpha \circ p_\alpha.
\end{equation}

For the index $\alpha \in A$ let us consider an index $\alpha '
\in A$, $\alpha ' \geqslant \alpha$, which satisfies the condition
of movability of the inverse system $\{X_\alpha , p_{\alpha \alpha
'}, A\}$. From (\ref{6}) we get
\begin{equation}
\label{7}
f=f_\alpha \circ p_{\alpha \alpha '} \circ p_{\alpha '}.
\end{equation}

Now let us prove that $X_{\alpha '}$, the homotopy classes
$p_\alpha :X \to X_{\alpha '}$ and $f_\alpha \circ p_{\alpha
\alpha '}:X_{\alpha '} \to Q$  satisfy condition $(*)$. Indeed,
let $Q''\in$ Ob($\mathcal{HCW}$) and $f'':X \to Q''$, $\eta ':Q''
\to Q$ homotopy classes, which satisfy the condition
\begin{equation}
\label{8}
f=\eta ' \circ f''.
\end{equation}

For the homotopy class $f'':X \to Q''$  there exist an index
$\alpha '' \in A$, $\alpha '' \geqslant \alpha$ and a homotopy
class $\tilde{f''}:X_{\alpha ''} \to Q''$   that
\begin{equation}
\label{9}
f''=\tilde{f''} \circ p_{\alpha ''}.
\end{equation}

It is clear that
\[
f_\alpha \circ p_{\alpha \alpha''} \circ p_{\alpha ''}=
\eta ' \circ \tilde{f''} \circ p_{\alpha ''}.
\]
Therefore, according to the definition 1 of "association", we can
find an index $\alpha ''' \in A$, $\alpha ''' \geqslant \alpha ''$
such that
\begin{equation}
\label{10}
f_\alpha \circ p_{\alpha \alpha''} \circ p_{\alpha '' \alpha '''}=
\eta ' \circ \tilde{f''} \circ p_{\alpha '' \alpha '''}.
\end{equation}

By the movability of the inverse system $\{X_\alpha , p_{\alpha
\alpha '}, A\}$, we can select the homotopy class $g:X_{\alpha '}
\to X_{\alpha '''}$ satisfying the condition
\begin{equation}\label{11}
p_{\alpha \alpha '}=p_{\alpha \alpha''} \circ p_{\alpha '' \alpha
'''} \circ g.
\end{equation}

\ \\[-6pt]

\begin{center}
\begin{picture}(100,100)
\put(50,50){$X$}
\put(0,50){$Q$}
\put(100,50){$X_{\alpha '''}$}
\put(25,0){$Q ''$}
\put(25,100){$X_{\alpha }$}
\put(75,0){$X_{\alpha ''}$}
\put(75,100){$X_{\alpha '}$}

\put(48,53){\vector(-1,0){38}}
\put(60,53){\vector(1,0){38}}
\put(52,58){\vector(-1,2){19}}
\put(58,47){\vector(1,-2){19}}
\put(52,47){\vector(-1,-2){19}}
\put(58,58){\vector(1,2){19}}
\put(25,97){\vector(-1,-2){19}}
\put(107,47){\vector(-1,-2){19}}
\put(87,97){\vector(1,-2){19}}
\put(25,7){\vector(-1,2){19}}
\put(75,103){\vector(-1,0){38}}
\put(75,3){\vector(-1,0){38}}

\put(2,75){$f_\alpha $} \put(32,75){$p_\alpha $}
\put(70,75){$p_{\alpha '}$} \put(103,75){$g$}

\put(2,25){$\eta '$} \put(32,25){$f ''$} \put(71,25){$p_{\alpha
''}$} \put(100,25){$p_{\alpha '' \alpha '''}$}

\put(25,57){$f$} \put(75,57){$p_{\alpha '''}$}
\put(50,107){$p_{\alpha \alpha '}$} \put(50,-12){$\tilde{f''}$}

\end{picture}

\end{center}

$$\text{Diagram 3.}$$

Let us define $\eta ''=\tilde{f''} \circ p_{\alpha '' \alpha '''}
\circ g$. It is remains to note that the homotopy class $\eta
'':X_{\alpha '} \to Q''$ satisfies the condition
\begin{equation}
\label{12} f_\alpha \circ p_{\alpha \alpha '}=\eta ' \circ \eta
''.
\end{equation}
Indeed:
\[
\eta ' \circ \eta ''=\eta ' \circ \tilde{f''} \circ p_{\alpha '' \alpha '''}
\circ g=f_\alpha \circ p_{\alpha \alpha ''} \circ p_{\alpha '' \alpha '''}
\circ g=f_\alpha \circ p_{\alpha \alpha '}.
\]

\end{proof}

\begin{remark}
The condition $(*)$ of Theorem 2 one can consider as a definition
of movability of topological space. \\ \\
\end{remark}

\end{document}